\documentclass[11pt]{amsart}

\usepackage{amsmath,amsfonts,amsthm,amssymb,amscd,stmaryrd,url,slashbox}
\usepackage{amssymb,slashbox}
\usepackage{multirow}

\newtheorem{theorem}{Theorem}[section]
\newtheorem{lemma}[theorem]{Lemma}

\theoremstyle{definition}

\theoremstyle{remark}

\newtheorem{cond}{Condition}

\numberwithin{equation}{section}

\newcommand{\myref}[1]{(\ref{#1})}

\newcommand{\C}{\mathbb{C}}
\renewcommand{\O}{{\mathcal{O}}}
\renewcommand{\a}{\alpha}
\renewcommand{\b}{\beta}
\renewcommand{\d }{\delta }

\newcommand{\g}{\gamma}
\renewcommand{\k}{\kappa}
\newcommand{\s}{\sigma}
\newcommand{\ep}{\epsilon}

\begin{document}

\title[Short intervals containing primes in arithmetic progressions]{Short effective intervals containing primes in arithmetic progressions and the seven cubes problem.}

\author{H. Kadiri}
\address{D\'epartement de Math\'ematiques et Statistique, Universit\'e de Montr\'eal, CP 6128 succ Centre-Ville, Montr\'eal, QC H3C 3J7, Canada}
 \email{kadiri@dms.umontreal.ca}

%\date{August 24, 2006 and, in revised form, .}


\begin{abstract}
For any $\ep>0$ and any non-exceptional modulus $q\ge3$, we prove that, for $x$ large enough ($x\ge \a_{\ep}\log^2 q$), the interval $\left[ e^x,e^{x+\ep }\right]$ contains a prime $p$ in any of the arithmetic progression modulo $q$.
We apply this result to establish that every integer $n$ larger than $\exp(71\,000)$ is a sum of seven cubes.
 \end{abstract}

\maketitle
%%%%%%%
\section{Introduction.}
Let $q\ge3$ be a non-exceptional modulus, $a$ a positive integer, $x>0$ and $\ep>0$ some real numbers. 
One way to establish that the interval $\left[ e^x,e^{x+\ep }\right] $ contains a prime $p \equiv a \pmod q$ would be to determine a condition on $x$ such that
\begin{equation}\label{theta}
\theta\left(e^{x+\ep };q,a\right)-\theta\left(e^x;q,a\right) = \sum_{\begin{substack}
{e^x< p \le e^{x+\ep }\\ p\equiv a [q]} \end{substack}}\log p
\end{equation}
is positive.
It will be convenient to work with Von Mangoldt's function
\begin{equation*}\Lambda(n) = \begin{cases} \log p & \hbox{ if } n= p^k, p \hbox{ prime},\\  0 & \hbox{ otherwise.} \end{cases} \end{equation*}
Showing that \myref{theta} is positive follows from showing that
\begin{equation*}
\psi \left( e^{x+\ep };q,a \right)- \psi \left( e^x;q,a \right) = \sum_{\begin{substack}
{e^x< n \le e^{x+\ep } \\ n\equiv a [q]} \end{substack}} \Lambda(n)
\end{equation*}
is larger than a positive constant times the error term between $\psi$ and $\theta$.
In \cite{MC}, following Rosser's method for $\psi(x)$ in \cite{Ros}, McCurley
approximated $\psi \left( e^x;q,a \right)$ via succesive integral
averaging. In fact, their methods amounts to weighting the primes
with a smooth function. Our approach will be to introduce directly a smooth positive weight into the difference $\psi \left( e^{x+\ep };q,a \right)- \psi \left( e^x;q,a \right)$:
\begin{equation}\label{sum0}
 \sum_{\begin{substack}{e^x< n \le e^{x+\ep } \\ n\equiv a [q]}\end{substack}} \dfrac{\Lambda(n)}{n}f(\log n).
\end{equation}
We choose the function $f$ so that it has compact support
contained in \(\left[x,x+\ep\right]\) and so that the peak of the function is near the prime we want to locate.
We have an explicit formula for the sum \myref{sum0}
\begin{equation}\label{tildepsi}
(1+o(1))\,\frac{F(0)}{\phi(q)} - \frac1{\phi(q)} \sum_{\chi \bmod q}
\overline{\chi}(a)\sum_{\varrho \in Z(\chi)} F(1-\varrho) ,
\end{equation}
where $F$ is the Laplace transform of $f$ and $Z(\chi)$ the set of
non-trivial zeros of $L(s,\chi)$. Note that this formula
generalizes the classical formula for $\psi$ (see pp.$121$-$122$ of \cite{Dav}):
\begin{equation*}
\psi(x;q,a)
= \frac{x}{\phi(q)} - \frac1{\phi(q)} \sum_{\chi \bmod q}
\overline{\chi} (a) \sum_{|\g|<T} \frac{x^{\varrho}}{\varrho}
+
\O\left(\frac{x\log^2(qx)}{\phi(q) T}  + \frac{x
e^{-c_1\sqrt{\log x}}}{\phi(q)}   \right).
\end{equation*}
The second argument relies on finding the largest real part for the zeros of the $L$-functions modulo $q$, in particular in the case of the low lying zeros.
The key result is due to Liu and Wang \cite{LiuWang}. It asserts that the zeros $\varrho =\b+i\g$ with $|\g|\le H$, except at most four of them, verify:
\begin{equation*}\b\le 1-\frac1{R_1\log (qH)},\ \hbox{ where } R_1=3.82.\end{equation*}
For the zeros of larger imaginary part, we use the latest effective result on the classical zero-free region (see \cite{Kad2}):
\begin{equation*}\b\le 1-\frac1{R\log (q|\g|)},\  \hbox{ where } R=6.50.\end{equation*}
We shall study the expression in \myref{tildepsi} with $x = \a \log^2q$. We deduce a lower bound for non-exceptional moduli $q$:
\begin{equation*}
\sum_{\begin{substack}{e^x< n \le e^{x+\ep } \\ n\equiv a
[q]}\end{substack}} \dfrac{\Lambda(n)}{n} \dfrac{f(\log
n)}{\|f\|_1} \ge \frac1{q} - (1+o(1))\frac{(\log H)\log(q^2H)}{2\pi\ep}\ q^{-\frac{\a}{R_1}\frac{\log q}{\log (qH)}},
\end{equation*}
where $H$ depends essentially on $\ep$ like: $H \asymp_q \ep^{-1}$.
From this we shall deduce that the sum on the primes is positive when:
\begin{equation*}
\a \ge R_1 \frac{\log(qH) }{\log^2q}\log \left(\frac{q(\log H)\log(q^2H)}{2\pi\ep}\right)(1+o(1)),
\end{equation*}
which gives values for $\alpha$ approaching $R_1$ as $\ep $ decreases. Our main result is the following:
\begin{theorem}\label{mainthm}
Let $q \ge 3$ be a non-exceptional modulus and let $(a,q)=1$. For any $\ep>0$, there exists an $\a >0$ such that, if $x \ge \a \log^2q$, then the interval $\left[ e^x,e^{x+\ep }\right] $ contains a prime $p \equiv a \pmod q$. 
Table \ref{table1} gives the values of $\a$ for various choices of $\ep$ and $q\ge q_0$. 
\end{theorem}
In section \ref{Section4}, we describe the general algorithm to compute $\a$ as a function of $q$ and $\ep$.
\begin{footnotesize}
\begin{table}[htbp]
\caption{} \label{table1} 
\renewcommand\arraystretch{1.1}
\noindent
\begin{tabular}{|c|c|c|c|c|c|c|} \hline
\backslashbox{$q_0$}{$\ep$} & 0.0001 & 0.001  & 0.01  & 0.1  & 1  & 10 \\
\hline
$5\cdot10^4$ & 19.228 & 15.550 & 12.245 & 9.4357 & 6.9684 & 4.8430 \\
$10^{10}$ & 9.8356 & 8.5912 & 7.4255 & 6.3398 & 5.3418 & 4.4761 \\
$10^{15}$ & 7.6121 & 6.8799 & 6.1816 & 5.5174 & 4.8905 & 4.3256 \\
$10^{20}$ & 6.5919 & 6.0799 & 5.5864 & 5.1114 & 4.6565 & 4.2373 \\
$10^{25}$ & 6.0079 & 5.6164 & 5.2364 & 4.8678 & 4.5116 & 4.1783 \\
$10^{30}$ & 5.6298 & 5.3137 & 5.0053 & 4.7047 & 4.4123 & 4.1357 \\
$10^{35}$ & 5.3649 & 5.0102 & 4.8411 & 4.5875 & 4.3396 & 4.1032 \\
$10^{40}$ & 5.1688 & 4.9414 & 4.7181 & 4.4989 & 4.2839 & 4.0776 \\
$10^{45}$ & 5.0178 & 4.8185 & 4.6225 & 4.4295 & 4.2398 & 4.0567 \\
$10^{50}$ & 4.8979 & 4.7205 & 4.5459 & 4.3737 & 4.2039 & 4.0394 \\
$10^{55}$ & 4.8003 & 4.6407 & 4.4832 & 4.3276 & 4.1740 & 4.0247 \\
$10^{60}$ & 4.7192 & 4.5742 & 4.4308 & 4.2890 & 4.1488 & 4.0121 \\
$10^{65}$ & 4.6509 & 4.5179 & 4.3864 & 4.2562 & 4.1272 & 4.0011 \\
$10^{70}$ & 4.5924 & 4.4697 & 4.3482 & 4.2278 & 4.1084 & 3.9915 \\
$10^{75}$ & 4.5418 & 4.4280 & 4.3151 & 4.2031 & 4.0920 & 3.9829 \\
$10^{80}$ & 4.4976 & 4.3914 & 4.2860 & 4.1814 & 4.0774 & 3.9753 \\
$10^{85}$ & 4.4587 & 4.3591 & 4.2603 & 4.1621 & 4.0645 & 3.9684 \\
$10^{90}$ & 4.4240 & 4.3304 & 4.2373 & 4.1448 & 4.0528 & 3.9622 \\
$10^{95}$ & 4.3931 & 4.3046 & 4.2168 & 4.1293 & 4.0423 & 3.9565 \\
$10^{100}$ & 4.3652 & 4.2815 & 4.1982 & 4.1153 & 4.0328 & 3.9513 \\
\hline
\end{tabular}
\end{table}\end{footnotesize}
In comparison, for $q\ge10^{30}$ and $\ep=\ln3$, McCurley's bounds on $\psi(x;q,a)$ would give $\alpha=10.690$ (see Theorem $1.2$ of \cite{KSMC}). 
With our new smoothing function, this result may first be improved to $\alpha=10.562$.
and with the new zero-free region ($R=6.50$ instead of $R=9.65$), to $\alpha=7.281$.
Using the region with a finite number of zeros ($R_1=3.82$), we finally obtain $\a= 4.401$. \\
Note that an explicit bound for the size of the least prime $p \equiv a \pmod q$, namely $P(a,q)$,
follows immediately:
\begin{equation}\label{1}
P(a,q) \le e^{\a \log^2 q}.
\end{equation}
In \cite{Wag}, Wagstaff computes the size of $P(a,q)$ for all possible arithmetic progressions up to modulus $5 \cdot 10^4$. For this reason, the data presented in Table \ref{table1} begins with moduli $q_0$ greater than $5 \cdot 10^4$.
There exists a stronger result than \myref{1} and we refer the reader to the work of Heath-Brown on the subject. In \cite{HB2}, he proved:
\begin{equation*}P(a,q) \ll q^{5.5}.\end{equation*}
Unfortunately, this is only valid for asymptotically large $q$. Moreover, if the proof is made effective, it is likely that this result would be weaker than \myref{1} in the range we are considering.
Also it can be applied to solve some effective problems.
We give an example in the second part of the article for which we will apply Theorem \ref{mainthm} for $q\ge 10^{32}$.
\\ \ \\
It concerns Waring's problem for sums of seven cubes.
Landau proved in $1909$ that every sufficiently large integer may be represented as a sum of eight nonnegative cubes. His proof used results on the representation of integers as a sum of three squares.
In $1943$, Linnik used a theorem on the representation of integers by ternary quadratic forms and proved in \cite{L} that it was also true with seven cubes.
In $1939$, Dickson completed Landau's statement by showing that all integers, except $23$ and $239$, is a sum of 8 cubes. It is widely expected that every integer $\ge 455$ is a sum of seven cubes.\\
In $1951$, Watson simplified Linnik's proof in \cite{Wa} by using a lemma establishing some conditions on $n$ to be represented as sum of seven cubes. This lemma has recently been improved by Ramar\'e in \cite{R}.
The main condition consists of finding prime integers in arithmetic progression as small as possible.
For example, McCurley found  \( n\ge \exp(1\,077\,334)\) in \cite{MC} and Ramar\'e \( n\ge \exp(205\,000)\) in \cite{R}. These authors use Chebyshev's estimates for $\theta(x;q,a)$ that McCurley previously established in \cite{KSMC}.
We replace this argument with our result concerning small intervals containing a prime mod $q$.
And since this assertion is only proven for non-exceptional modulus, we give an explicit description of the scarcity of exceptional moduli. We prove in section \ref{7cubes}:
\begin{theorem}\label{cubes}
Every integer $n$ larger than $N_0=\exp(71\,000)$ is a sum of seven cubes.
\end{theorem}
%%%%%%%
%%%%%%%
\section{Preliminary lemmas.}\label{Section2}\label{preliminarylemmas}
%%%%%%%
%%%%%%%
\subsection{Zeros of the Dirichlet $L$-functions.}
The proof depends essentially on the distribution of the zeros of Dirichlet $L$-functions.
The first theorem states an explicit zero free region for all moduli $q$, even for those of not too large size.
\begin{theorem}(Theorem $1.1$ of \cite{Kad2})\label{th1}
Let $q$ be an integer, $q \ge 3$, and let $\mathcal{L}_q(s)$ be the product of Dirichlet $L$-functions
modulo $q$. Then $\mathcal{L}_q(s)$ has at most one zero in
the region :
\begin{equation*}
\sigma \ge 1-\frac1{R \log \left(q\max(1,|t|)\right)}, \hbox{ where }   R= 6.50.
\end{equation*}
Such a zero, if it exists, is real, simple and corresponds to a real non-principal character modulo $q$. We
shall refer to it as an exceptional zero and $q$ as an exceptional modulus.
\end{theorem}
The next theorem illustrates the fact that the zeros do not cluster near the one-line. In fact, there are few of them:
\begin{theorem}[Theorem $1$ of \cite{LiuWang}]\label{LiuWang}
Suppose $q$ is an integer satisfying $1\le q \le x$ and $x$ a real number, $x \ge 8\cdot 10^9$. Then the function $\mathcal{L}_q(s)$ has at most four zeros in
the region
\begin{equation*}
|\Im s| \le x/q,\quad \sigma \ge 1-\frac1{R_1 \log x}, \hbox{ where }  R_1= 3.82. 
\end{equation*}
\end{theorem}
We will apply this theorem for the case $x=q$ and $x=qH$.\\
We describe explicitly the following phenomenon: the exceptional zero tends to repel the zeros of close conductor.
\begin{theorem}[Theorem $1.3$ of \cite{Kad2}]
\label{repulsion}
If $\chi_1$ and $\chi_2$ are two distinct real primitive characters modulo $q_1$ and $q_2$ respectively and if $\beta_1$ and $\beta_2$ are real zeros of $L(s,\chi_1)$ and $L(s,\chi_2)$ respectively, then:
\begin{equation*}
\min\left(\beta_1,\beta_2 \right) \le 1-\frac1{ R_2  \log(q_1q_2)},\hbox{ where } R_2 = 2.05. 
\end{equation*}
\end{theorem}
When $q_1< q_2$, then both $q_1$ and $q_2$ cannot be exceptional, unless
$
q_2 \ge q_1^{2.12}.
$
The next theorems gives explicit density for the zeros associated to each character $\chi$ modulo $q$ (see p.$267$ of \cite{KSMC}).
\begin{lemma}\label{McCurley}
Let $T\ge1$. We denote $N(T,\chi)$ the number of zeros of the Dirichlet
$L$-function $L(s,\chi)$ in the rectangle $ \displaystyle{ \{s\in \C: 0 \le \Re s \le 1 , | \Im s | \le T \}}$. \\
Then $ N(T,\chi) = P(T) + r(T),$ 
with\begin{equation*}P(T):= \frac{T}{\pi}\log \frac{qT}{2\pi e},\ |r(T)| \le R(T):= a_1 \log(qT) + a_2, \ a_1=0.92 ,\ a_2 = 5.37.\end{equation*}
\end{lemma}
The next lemma establish a bound for 
\begin{equation*}\mathcal{S}(H,\chi):=\sum_{\begin{substack}{1<|\g|<H\\ L(\b+i\g,\chi)=0}\end{substack}} \frac{1}{|\g|} .\end{equation*}
\begin{lemma}\label{sum1overRho}
Let $q$ be the conductor of $\chi$ and $H\ge1$, then
$\mathcal{S}(H,\chi) \le \tilde E(H)$ , where
\begin{align*}
\tilde E(H):=& \frac{1}{\pi} (\log q)(\log H)
+ \frac1{2\pi} \log^2H
+ \left(\frac{1}{\pi} + a_1\right) \log q
- \frac{\log(2\pi)}{\pi} \log H
\\&- \frac{1}{\pi}\log (2\pi e) + a_2 + a_1
- \frac{a_1}H.
\end{align*}
\end{lemma}
\begin{proof}
We have:
\begin{equation*}
\mathcal{S}(H,\chi)
=  \frac{N(H,\chi)}{H} - N(1,\chi) + \int_1^{H}\frac{N(t,\chi)}{t^2} d\,t .
\end{equation*}
We use Lemma \ref{McCurley} to bound $N(t,\chi)$ in the integral and we integrate by parts to obtain:
\begin{equation*}
\mathcal{S}(H,\chi)
\le P(1)+R(1) + \int_1^{H} \frac{P'(t)+R'(t)}{t} d\,t.
\end{equation*}
We conclude by computing the last integral:
\begin{equation*}
\int_1^{H} \frac{P'(t)+R'(t)}{t} d\,t 
= 
\frac{1}{\pi} (\log q)(\log H)
+ \frac1{2\pi} \log^2H
- \frac{\log(2\pi)}{\pi} \log H
+ a_1- \frac{a_1}H.
\end{equation*}
\end{proof}
%%%%%%%
%%%%%%%
\subsection{Bounds for $\Gamma'/\Gamma(s)$.}
\begin{lemma}
\label{lem6}
If $T\ge0$, then
\begin{equation*}
\left|\Re\dfrac{\Gamma'}{\Gamma}\left(\frac{2-\chi(-1)}{4}+i\frac{T}{2}\right)\right|\\
\le U(T) := \log\left( 6(T+12)\right) .
\end{equation*}
\end{lemma}
\begin{proof}See \cite{Kad2}.\end{proof}
%%%%%%%
%%%%%%%
\subsection{Properties of the weight function.}
Let $m$ be a positive integer, $L$ and $\ep $ some positive
constants. Our choice for $f$ is inspired by the function Ramar\'e
and Saouter used in p.$17$ of \cite{RS}. They call $f$
$m-$admissible when it satisfies:
\begin{itemize}
\item $f$ is an $m-$times differentiable function,
\item $f^{(k)}(0) = f^{(k)}(1) =0$ if $0 \le k \le m-1$,
\item $f \ge0$,
\item $f$ is non-identically zero.
\end{itemize}
The specific function we use is
\begin{equation}\label{fdef}
f(t) = (t-L)^m\, (L+\ep -t)^m, \hbox{ if } L\le  t \le L+\ep ,
\end{equation}
and \(f(t)=0\) otherwise. Furthermore, we notice that $f$ and
its derivative are symmetric with respect to $L+\ep /2$.
\begin{lemma}\label{lem1}
\begin{align}
& \label{norm1mf}
\frac{\|f^{(m)}\|_2}{\|f\|_1}=  \frac{\mu_m}{\ep ^{m+1/2}},
\hbox{ with }
\mu_m = \frac{(2m+1)!}{m!\sqrt{2m+1}}.
\\ & \label{norminf}
\frac{\|f\|_{\infty}}{\|f\|_1}= \frac{\nu_m}{\ep },
\hbox{ with }
\nu_m = \frac{(2m+1)!}{4^m (m!)^2}.
\end{align}
\end{lemma}
\begin{proof}
This is an exercise, or else, see p.$17$ of \cite{RS}, for the derivation.
\end{proof}
Let $F$ be the Laplace transform of the function $f$ as defined in \myref{fdef}:
\begin{equation*}
F(s)= \int_{L}^{L+\ep } f(t) e^{-st} d\,t.
\end{equation*}
\begin{lemma}\label{lem2}
If $\s \ge 0$:
\begin{align}
& \label{Laplacef+2}
F(\s) \ge e^{ -\s (L+\ep)}\|f\|_1 ,\\
& \label{Laplacef+2bis} 
|F(s)| \le e^{ -\s L}\|f\|_1,
\\ & \label{Laplacef1+}
|F(s)| 
\le \frac{e^{-\s L}}{|s|}  \frac{2(2m+1)}{\ep m } \|f\|_1.
\\ & \label{Laplacefm+}
|F(s)| \le \sqrt{\ep }\, e^{ -\s L} \frac{\|f^{(m)}\|_2 }{|s|^m}.
\end{align}
\end{lemma}
\begin{proof}
The proof  makes use of the symmetry of $f$ and $f^{(m)}$ and, in the last case, the Cauchy-Schwarz inequality.
\end{proof}
%%%%%%%
%%%%%%%
\subsection{An explicit formula.}
Let $q$ be an integer, $q\ge3$. For each character $\chi$ modulo $q$, we denote $\chi_1$ the primitive character asociated to $\chi$.
\begin{lemma}\label{S11}
Let $f$ be the function given by \myref{fdef}.
Then
\begin{equation*}
\begin{split}
\sum_{\chi \bmod q} \overline{\chi(a)} \sum_{n\ge1}
\frac{\Lambda(n) \chi_1(n)}{n} f(\log n) =& F(0) + c(a,q) F(1)
\\& - \sum_{\chi \bmod q} \overline{\chi(a)} \sum_{\varrho \in Z(\chi_1)} F(1-\varrho) + I(a,q) .
\end{split}
\end{equation*}
where $c(a,q) \ge \frac12$,
$ Z(\chi_1)$ is the set of non-trivial zeros of $L(s,\chi_1)$ and\\
$ \displaystyle{ I(a,q) = \frac1{2\pi} \int_{-\infty}^{+\infty}
\sum_{\chi \bmod q} \overline{\chi(a)} \Re \frac{\Gamma'}{\Gamma}
\left(\frac{2-\chi_1(-1)}{4}+i\frac{T}2 \right) F\left(1/2-iT\right)
d\,T.}$
\end{lemma}
\begin{proof}
This is a special case of the explicit formula of Theorem 3.1,
p.$314$ of \cite{Kad1}, applied to the smooth function
\( \phi (x) = f(x) e^{- x} \) if \(x \ge 0\) and \( \phi (x) = 0 \) otherwise.
The constant $c(a,q)$ is given by
\begin{equation*}
 \frac14 \sum_{\chi \bmod q } \overline{\chi(a)} + \frac14 \sum_{\chi \bmod q } \overline{\chi(a)}\chi_1(-1) + \frac12 \ge \frac12.
\end{equation*}
(see p.$414$ of \cite{RR} for details).
\end{proof}
%%%%%%%
%%%%%%%
\section{Main lemma.}\label{proof}
We place our study in the case of modulus $q$ not studied by Wagstaff, that is to say for those larger than $5 \cdot 10^4$.
Let $\ep>0$, $H\ge1$ and $\a$ be positive reals such that
\begin{equation}\label{Condalf}
 \alpha < R\left(\frac{\log(qH)}{\log q}\right)^2 .
\end{equation}
We define $L$ as a parameter depending only on $q$:
\begin{equation*}
L:=\alpha \log^2 q .
\end{equation*}
Throughout the paper, $\varrho=\b+i\g$ always stands for a non-trivial zero of a Dirichlet $L$-function.
%%%%%%%
We prove in this section that
\begin{lemma}\label{mainlem}
If $m\ge3$, $\ep>0$ and if $\a$ is satisfying the condition \myref{Condalf}, then the sum over the primes
$ \displaystyle{
 \Sigma(a,q)  := \sum_{{\begin{substack} {p \equiv a \bmod q}  \end{substack}}}\frac{\log p}{p} f(\log p) }$
satisfies
\begin{equation*}
\frac{\Sigma(a,q)}{\|f\|_1}\ge \frac1q - r(\a,\ep,H,m,q),
\end{equation*}
where $ \displaystyle{r := \sum_{i=1}^5r_i}$ and the $r_i$'s are given by \myref{ep1} , \myref{ep2} , \myref{ep3} , \myref{ep4} and \myref{ep5}.
\end{lemma}
Note that $\Sigma(a,q)$ is actually close to the sum appearing in
Lemma \ref{S11}:
\begin{equation*}
 \Sigma = \Sigma_{11} + \Sigma_{12} - \Sigma_2
\end{equation*}
\begin{align*}
\hbox{with }\qquad
& \Sigma_{11}(a,q) := \frac{1}{\phi(q)}
\sum_{\chi \bmod q} \overline{\chi(a)} \sum_{n\ge1}
\frac{\Lambda(n) \chi_1(n)}{n} f(\log n),
\\
& \Sigma_{12}(a,q) := \frac{1}{\phi(q)} \sum_{\chi \bmod q} \overline{\chi(a)} \sum_{n\ge1} \frac{\Lambda(n) (\chi(n)-\chi_1(n))}{n} f(\log n),
\end{align*}
where $\chi_1$ is the primitive character associated
to $\chi$, and
\begin{equation*}
 \Sigma_2(a,q) :=  \sum_{{\begin{substack} {k\ge 2 \\ p^k \equiv a \bmod q}  \end{substack}}}\frac{\log p}{p^{k}} f(k\log
 p).
\end{equation*}
We will prove in sections \ref{sigma12} and \ref{sigma2} that the
two last sums are small in comparison to $\Sigma_{11}$.
We use Lemma \ref{S11} to bound $\Sigma_{11}$:
\begin{equation*}
\Sigma(a,q) \ge A_1(q) - A_2^-(q)- A_2^+(q) - A_3(q) -\left|\Sigma_{12}(a,q)\right|-\left|\Sigma_{2}(a,q)\right|,
\end{equation*}
\begin{align*}
\hbox{with  }
&A_1(q) := \frac1{\phi (q)}\left(F(0) + c(a,q) F(1)\right) , \\
&A_2^-(q) := \frac1{\phi (q)}\sum_{\chi \bmod q}
 \sideset{}{^{\prime}}\sum_{{\begin{substack} {\varrho \in Z(\chi) \\ |\g|\le H}  \end{substack}}} \left| F(1-\varrho) + F(\overline{\varrho})\right|,\\
&A_2^+(q) := \frac1{\phi (q)}\sum_{\chi \bmod q}
 \sideset{}{^{\prime}}\sum_{{\begin{substack} {\varrho \in Z(\chi) \\ |\g| > H}  \end{substack}}} \left| F(1-\varrho) + F(\overline{\varrho}) \right|,\\
&A_3(q) := \frac1{\phi (q)}|I(a,q)|,
\end{align*}
where $ \displaystyle{ \sideset{}{^{\prime}}\sum_{\beta } =
\sum_{ 1/2<\b <1} + \frac1{2} \sum_{\b =1/2} }$ (we use the symmetry property of the zeros).
We extend the sum over the zeros of $L(s,\chi_1)$ to the zeros of $L(s,\chi)$ to simplify our argument.\\
The sections \ref{Section31} to \ref{Section33} study these $A_i$'s.
%%%%%%%
\subsection{Study of $A_1$.}\label{Section31}
It is immediate that:
\begin{equation}
\label{boundA1fin} \dfrac{A_1(q)}{\|f\|_1}\ge \frac1q.
\end{equation}
%%%%%%%
%%%%%%%
\subsection{Study of $A_2^-$.} \label{Section32}
Since we are in the case where $q$ is non-exceptional, we do not
worry about the existence of a Siegel zero. Thanks to theorems
\ref{th1} and \ref{LiuWang} we can split the sum $A_2^{-}$ as
follow
\begin{align*}
A_2^-(q) =& \frac1{\phi(q)} \sum_{k=1}^8 \left| F(1-\varrho_k) +
F(\overline{\varrho}_k) \right| 
+ \frac1{\phi (q)}\sum_{\chi \bmod q}
\sideset{}{^{\prime}}\sum_{\begin{substack} { |\g | \le 1 \\ \b \le 1-\frac1{R_1\log q} }\end{substack}}  \left|
F(1-\varrho) + F(\overline{\varrho}) \right|
\\ &+ \frac1{\phi (q)}\sum_{\chi \bmod q}
\sideset{}{^{\prime}}\sum_{\begin{substack} { 1< |\g | \le H \\ \b \le 1-\frac1{R_1\log (qH)} }\end{substack}}  \left|
F(1-\varrho) + F(\overline{\varrho}) \right|,
\end{align*}
where the zeros $ \varrho_k=\b_k+i\g_k$ satisfy:
\begin{align*}
& |\g_{k}| \le 1 \hbox{ and } 1-\frac1{R_1\log q} \le \b_k \le 1-\frac1{R\log q}\ , \hbox{ for } k=1,2,3,4 ,\\
& |\g_{k}| \le H \hbox{ and } 1-\frac1{R_1\log (qH)} \le \b_k \le 1-\frac1{R\log (qH)}\ , \hbox{ for } k=5,6,7,8.
\end{align*}
%%%
We use the inequalities \myref{Laplacef+2bis} and \myref{Laplacef1+} for the first and second line respectively:
\begin{align*}
& \left| F(1-\varrho) + F(\overline{\varrho}) \right| 
\le \left( e^{ -\b L} + e^{ -(1-\b) L}\right)\|f\|_1\\
\hbox{and}& \left| F(1-\varrho) + F(\overline{\varrho}) \right| \le \frac{1}{|\g|}\left( e^{ -\b L} + e^{ -(1-\b) L} \right) \frac{2(2m+1)}{\ep m }\|f\|_1.
\end{align*}
\begin{align*}
\hbox{Then  }\quad
A_2^-(q) 
\le & \frac{4}{\phi(q)} b_2(\alpha,R,q) \left( q^{ -\frac{\a}{R}} + q^{ -\frac{\a}{R}\frac{\log q}{\log(qH)}} \right) \|f\|_1  
\\&  + \frac{1}2 \,b_2(\alpha,R_1,q) q^{ -\frac{\a}{R_1}} \|f\|_1 \,\max_{\chi \bmod q}N(1,\chi)
\\& + b_2(\alpha,R_1,q)  \frac{2m+1}{\ep m } \|f\|_1 \, q^{ -\frac{\a}{R_1}\frac{\log q}{\log (qH)}} \, \max_{\chi \bmod q} \left( \sideset{}{^{\prime}}\sum_{\begin{substack} {1< |\g | \le H }\end{substack}} \frac{1}{|\g|}\right) ,
\end{align*}
where   $b_2(\alpha,r,q) := 1+  q^{-\a \log q+\frac{2\a}{r} }$.
We conclude by bounding the sum over the zeros like in Lemma \ref{sum1overRho}, $N(1,\chi)$ like in Lemma \ref{McCurley} and $\phi(q)$ like in page $72$ of \cite{RS1}: 
\begin{equation*}\frac{q}{\phi(q)} < e^{C}\log\log q + \frac{2.51}{\log\log q} \hbox{  for  } q\ge3, \end{equation*}
where $C$ stands for the Euler constant. 
 \begin{align*}
\frac{A_2^-(q)}{\|f\|_1}
\le & 4 b_2(\alpha,R,q) \left( q^{ -1-\frac{\a}{R}} + q^{ -1-\frac{\a}{R}\frac{\log q}{\log(qH)}} \right)\left(e^{C}\log\log q + \frac{2.51}{\log\log q}\right)
\\& + \Big(\frac{1+a_1\pi}{2\pi}\log q - \frac{\log(2\pi e)+ a_2\pi}{2\pi}  \Big) b_2(\alpha,R_1,q) q^{ -\frac{\a}{R_1}} 
\\& +  \frac{(2m+1) \tilde E(H)}{2\ep m }b_2(\alpha,R_1,q)q^{ -\frac{\a}{R_1}\frac{\log q}{\log (qH)}}.
\end{align*}
We obtain
\begin{equation} \label{boundA2-fin} 
\dfrac{A_2^-(q)}{\|f\|_1} \le r_1(\a,\ep,H,m,q) ,
\end{equation}
\begin{equation}\label{ep1}
\begin{split}
& r_1(\a,\ep,H,m,q) := 
b_2(\alpha,R_1,q)  \frac{(2m+1) }{2\ep m }\Bigg(
\frac{(\log H)(\log (q^2H))}{2\pi}
\\&+ \Big(\frac{1}{\pi} + a_1\Big) \log q
- \frac{\log(2\pi)}{\pi} \log H
- \frac{\log (2\pi e)}{\pi} + a_2 + a_1
- \frac{a_1}H\Bigg) q^{ -\frac{\a}{R_1}\frac{\log q}{\log (qH)}}
\\& +  b_2(\alpha,R_1,q)  \Big(\frac{1+a_1\pi}{2\pi}\log q - \frac{\log(2\pi e)+ a_2\pi}{2\pi}  \Big)q^{ -\frac{\a}{R_1}} 
\\& + 4 b_2(\alpha,R,q) \left(e^{C}\log\log q + \frac{2.51}{\log\log q}\right) \left( q^{ -1-\frac{\a}{R}} + q^{ -1-\frac{\a}{R}\frac{\log q}{\log(qH)}} \right).
\end{split}
\end{equation}
%%%%%%%
%%%%%%%
\subsection{Study of $A_2^+$.}\label{Section32bis}
Theorem \ref{th1} allows us to restrict
$\sideset{}{^{\prime}}\sum$ to the zeros in the region:
\begin{equation*}
|\g| > H,\  1/2 \le \b \le 1- \frac1{R\log (q|\g|)}.
\end{equation*}
We use \myref{Laplacefm+} to bound $|F(1-\varrho)|$ and$|F(\overline{\varrho})|$ 
\begin{multline*}
\left| F(1-\varrho) + F(\overline{\varrho}) \right|
\\ \le \sqrt{\ep } \|f^{(m)}\|_2
\Big[ \exp\Big( \frac{-L}{R\log(q|\g|)} \Big) +\exp\Big( -L \Big( 1-\frac1{R\log (qH)} \Big) \Big) \Big]\frac{1}{|\g|^m}.
\end{multline*}
We follow Lemma 4.1.3. and Lemma 4.2.1. of \cite{RR} and obtain that if $L\le R\log^2(qH)$, then:
\begin{equation*}
\sideset{}{^{\prime}}\sum_{\begin{substack} {\varrho \in Z(\chi_1)
\\ |\g | > H} \end{substack}} \frac{\exp\Big( \frac{-L}{R\log(q|\g|)} \Big)}{|\g|^m}
\le \frac{\tilde A + \tilde B}{2}
\quad \hbox{   and}\quad 
\sideset{}{^{\prime}}\sum_{\begin{substack} {\varrho \in Z(\chi_1)
\\ |\g | > H} \end{substack}} \frac{1}{|\g|^m} \le
\frac{\tilde C+\tilde D}{2}, 
\end{equation*}
with
\begin{align*}
&\tilde A := 
\frac{1}{\pi(m-2)H^{m-1}} \exp\left(-\frac{L}{R\log(qH)}\right) \left(\log\frac{qH}{2\pi}+\frac1{m-2}+\frac{a_1}{(m-1)H}\right),\\
&\tilde B := \frac{2(a_1\log(qH)+a_2)}{H^m}\exp\left(-\frac{L}{R\log(qH)}\right),\\
&\tilde C := \frac1{\pi(m-1)H^{m-1}}\left(\log\frac{qH}{2\pi}+\frac1{m-1}\right) ,\\
&\tilde D := \frac{2a_1\log(qH)+2a_2+\frac{a_1}{m}}{H^m}.
\end{align*}
We deduce the bound:
\begin{equation}\label{sumA2+}
\frac{A_2^+(q)}{\sqrt{\ep } \|f^{(m)}\|_2} \le \frac{\tilde A + \tilde B}{2} + \frac{\tilde C + \tilde D}{2} \exp\Big( -L \Big( 1-\frac1{R\log (qH)} \Big) \Big) 
\end{equation}
and together with \myref{norm1mf}:
\begin{equation}
\label{boundA2+fin} \dfrac{A_2^+(q)}{\|f\|_1} \le r_2(\a,\ep,H,m,q ) 
\end{equation}
\begin{equation}\label{ep2}
\begin{split}
&r_2(\a,\ep,H,m,q ):=
 q^{-\frac{\a}{R} \frac{\log q}{\log(qH)}} \frac{\mu_m}{(H\ep)^m} 
\Big[ 
\frac{H\log\frac{qH}{2\pi}}{2\pi(m-2)}
+\frac{H}{2\pi(m-2)^2}
\\&+\frac{a_1}{2\pi(m-2)(m-1)}
+ a_1\log(qH)+ a_2 \Big]
+ q^{-\a \log q + \frac{\a}R \frac{\log q}{\log(qH)}} \frac{\mu_m}{(H\ep)^m} 
\\&\Big[ \frac{H}{2\pi(m-1)}\left(\log\frac{qH}{2\pi}+\frac1{m-1}\right) 
+ a_1\log(qH)+a_2+\frac{a_1}{2m} \Big] . 
\end{split}
\end{equation}
%%%%%%%
%%%%%%%
\subsection{Study of $A_3$}\label{Section33}
\begin{equation*}
A_3(q) \le  \frac1{2\pi} \int_{-\infty}^{+\infty} \left| \Re
\dfrac{\Gamma'}{\Gamma} \left(\frac{2-\chi_1(-1)}{4}+i\frac{T}2
\right)\right| \left|F\left(1/2-iT\right)\right| d\,T.
\end{equation*}
We use Lemma \ref{lem6} to bound $\Gamma'/\Gamma$,
 \myref{Laplacef+2bis} to bound $F$ when $T\le 1$ and
\myref{Laplacefm+} otherwise:
\begin{align}
&\label{boundA3fin}
\frac{A_3(q)}{\|f\|_1} \le r_3(\a,\ep,m,q ),
\\&\label{ep3}
 r_3(\a,\ep,m,q ) := \frac1{2\pi}\left(  J_0 + \frac{\mu_m J(m)}{\ep ^m }  \right)q^{-\frac{\alpha}2 \log q}
\end{align}
with $\displaystyle{ J_0 := \int_{|T|\le1} \log(6(|T|+12)) d\,T }$  and  $\displaystyle{ J(m) := \int_{|T|>1} \frac{\log(6(|T|+12))}{|T|^m}d\,T}$.
%%%%%%%
%%%%%%%
\subsection{Study of $\Sigma_{12}(a,q)$.}\label{sigma12}
For $n$ fixed, we denote $Q_n$ the largest divisor of $q$ coprime with $n$. Then
\begin{equation*}\frac1{\phi(q)}\sum_{\chi \bmod q} \chi_1(n) \overline{\chi(a)}
=\begin{cases}
\frac{\phi(Q_n)}{\phi(q)} &\hbox{if }\ n \equiv a\ (\bmod\ Q_n),\\
0 &\hbox{else }.
\end{cases}\end{equation*}
For a proof, see p$414$ of \cite{RR}.
It implies that
\begin{equation*}
\Sigma_{12}(a,q)
= \sum_{\begin{substack} {n \equiv a \bmod Q_n \\ Q_n <q} \end{substack}} \frac{\phi(Q_n)}{\phi(q)} \frac{\Lambda(n)f(\log n)}{n}.
\end{equation*}
In this sum, we have
\begin{equation*}
\frac{\phi(Q_n)}{\phi(q)} 
= \frac1{p^{\nu_{p}(q)-1}(p-1)}
\end{equation*}
since $n$ is a prime power, $n=p^k$, coprime with $Q_n$ but not with $q$. Therefore
\begin{equation}\label{X1}
\Sigma_{12}(a,q)
\le \|f\|_{\infty} \sum_{p^{\nu_{p}(q)}| q} \frac{\log p}{p^{\nu_{p}(q)-1}(p-1)} \sum_{e^L < p^k < e^{L+\ep}} \frac{1}{p^k}.
\end{equation}
We compute the geometric sum
\begin{equation}\label{X2}
\sum_{e^L < p^k < e^{L+\ep}} \frac{1}{p^k}
\le \sum_{k \ge \left[\frac{L}{\log p}\right]+1}\frac{1}{p^k}
= \frac{e^{-L}}{p-1}.
\end{equation}
We reinsert the last bound in the summand and split the obtained sum:
\begin{equation}\label{X3}
\sum_{p^{\nu_{p}(q)}| q} \frac{\log p}{p^{\nu_{p}(q)-1}(p-1)^2}
\le \sum_{p| q} \frac{\log p}{(p-1)^2}+\sum_{p^j| q,j\ge2} \frac{\log p}{p^{j-1}(p-1)^2},
\end{equation}
where
\begin{equation}\label{X4}
\sum_{p^j| q,j\ge2} \frac{\log p}{p^{j-1}(p-1)^2}
\le
\sum_{p| q} \frac{\log p}{(p-1)^2} \sum_{j\ge2} \frac{1}{p^{j-1}}
= \sum_{p| q} \frac{\log p}{(p-1)^3}
\end{equation}
and
\begin{equation}\label{X5}
\sum_{p\ge2} \log p \left( \frac{1}{(p-1)^2} + \frac{1}{(p-1)^3} \right)
\le 2.10.
\end{equation}
Together with \myref{norminf} and \myref{X1} to \myref{X5}, we conclude that
\begin{align}
& \label{boundS12-c6}
\frac{\left|\Sigma_{12}(a,q)\right|}{ \|f\|_{1}}
\le 2.10\frac{\|f\|_{\infty}}{\|f\|_{1}} e^{-L} \le r_4(\a,\ep,m,q) ,\\&\label{ep4}
r_4(\a,\ep,m,q) := 2.10\, \frac{\nu_m}{\ep }\, q^{-\a \log q}.
\end{align}
%%%%%%%%%
%%%%%%%
\subsection{Study of $\Sigma_{2}(a,q)$.}\label{sigma2}
\begin{equation*}
\left|\Sigma_2(a,q)\right| 
= \sum_{{\begin{substack} {k\ge 2 \\ p^k \equiv a \bmod q}  \end{substack}}}\frac{\log p}{p^{k}} f(k\log p)
 \le \|f\|_{\infty} \sum_{2 \le p \le e^{\frac{L+\ep }2}}\log p\ \sum_{e^L < p^k < e^{L+\ep}} \frac{1}{p^k}
\end{equation*}
From \myref{X2} and
\begin{equation*}\sum_{2 \le p \le e^{\frac{L+\ep }2}}\frac{\log p}{p-1} 
\le 2\log\left(e^{\frac{L+\ep }2}\right) = L+\ep \hbox{ (see equation $(3.24)$ of \cite{RS1})},\end{equation*}
it follows:
\begin{align}
&\label{boundS2-c7}
\dfrac{\left|\Sigma_2(a,q)\right|}{ \|f\|_{1}}
\le \frac{\|f\|_{\infty}}{\|f\|_{1}} (L+\ep) e^{-L}
\le r_5(\a,\ep,m,q),
\\&\label{ep5}
r_5(\a,\ep,m,q) := \frac{\nu_m}{\ep } (\a \log^2q + \ep) q^{-\a \log q}.
\end{align}
%%%%%%%
%%%%%%%
\section{Proof of Theorem \ref{mainthm}}\label{Section4}
We gather the inequalities \myref{boundA1fin}, \myref{boundA2-fin}, \myref{boundA2+fin}, \myref{boundA3fin}, \myref{boundS12-c6} and \myref{boundS2-c7} and obtain:
\begin{equation*}
\frac{\Sigma(a,q)}{\|f\|_1}\ge  q^{-1} - r(\a,\ep,H,m,q ) \ge  q^{-1}\left( 1 - q_0 r(\a,\ep,H,m,q_0 ) \right).
\end{equation*}
Let $u\in [0.001,0.2]$, $q\ge q_0$ with $q_0=5\cdot10^4, 10^{10},...,10^{100}$ and $\ep=10^{-3},10^{-2},...,10$, be fixed.
We will choose $H$ and $m$ such that $\a$ is as small as possible and satisfies
\begin{equation}
\label{Cond2}
1 - q_0r(\a,\ep,H,m,q_0 )= 10^{-6}
\end{equation}
and $r_1$ and $r_2$ are of comparable size:
\begin{equation}
\label{Cond1}
r_2(\a,\ep,H,m,q_0 ) =ur_1(\a,\ep,H,m,q_0 ).
\end{equation}
We approximate $r_{1}$, $r_{2}$ and $r$ with $\tilde r_{1}$, $\tilde r_{2}$ and $\tilde r_{1} + \tilde r_{2}=(1+u)\tilde r_{1}$ respectively, where
\begin{align*}
& \tilde r_1(\a,H,m) :=  \frac{(2m+1) (\log H)(\log (q_0^2H))}{4\pi\ep m} q_0^{ -\frac{\a}{R_1}\frac{\log q_0}{\log (q_0H)}}
,\\
& \tilde r_2(\a,H,m)
:= \frac{H\log (q_0H)}{\pi\sqrt{m}} q_0^{-\frac{\a}{R} \frac{\log q_0}{\log(q_0H)}}\left(\frac{4m}{eH\epsilon}\right)^m .
\end{align*}
We approximate \myref{Cond2} by the equation
\begin{equation*}
1 - q_0(1+u)\tilde r_1(\a,H,m)= 10^{-6}.
\end{equation*}
Its solution is close to
\begin{equation}\label{tildealfa}
\tilde \a(H,m) 
 := R_1 \frac{\log (q_0H)}{\log^2 q_0} \log\left( \frac{ q_0 (\log H)(\log (q_0^2H))}{2 \pi\ep} \right) .
\end{equation}
It remains to find appropriate values of $H$ which will satisfy \myref{Cond1}.
The solution of the equation
\begin{equation*}
\tilde r_2(\tilde \a(H,m),H,m) = u \tilde r_1(\tilde \a(H,m),H,m)
\end{equation*}
is close to
\begin{equation} \label{tildeH}
\tilde H (m) 
=\frac1{\ep} \left(\frac{4 }{u \sqrt{m}} \left(\frac{4m}{e}\right)^m \left( \frac{ q_0 \log q_0}{4 \pi\ep } \right)^{1-\frac{R_1}{R}}  \right)^{\frac1{m-1}}.
 \end{equation}
 We minimize the value of $\tilde \a(\tilde H(m),m)$ and find that $m$ is close to
\begin{equation*}
\tilde m:= \frac{1}{2} + \log \left(\frac{16 }{u}\left( \frac{ q_0 \log q_0}{4 \pi\ep } \right)^{1-\frac{R_1}{R}} \right).
\end{equation*}
We now describe the algorithm to compute $\a$.
For $u$ and $m$ fixed (the value of $m$ is chosen close to $\tilde m$),
\begin{itemize}
\item We compute $\tilde H(m)$ and $\tilde \a (\tilde H(m),m)$ as given in \myref{tildeH} and \myref{tildealfa} respectively.
\item We choose for $H$ the value of the solution of the following approximation of equation \myref{Cond1}:
\begin{equation*}
r_2(\tilde \a (\tilde H(m),m),H,m )=ur_1(\tilde \a (\tilde H(m),m),H,m ).
\end{equation*}
With this value for $H$, we solve \myref{Cond2} with respect to $\a$. It is not difficult to see that the function $r(\a,\ep,H,m,q )$ decreases when $\a$ increases. Therefore we are insured of the uniqueness of the solution of the equation.
\item We choose $u$ and $m$ so that the value of $\a$ is as small as possible.
\end{itemize}
Table \ref{computations} records the values of the parameters $m$, $H$ and $u$. They have been rounded up in the last decimal place.
\\
For the next section, we will use the following result:
when $q\ge10^{32}$, $\ep=1.9$, then $u=0.022$, $H=80.8$, $m=38$ and $\a=4.3060$.
%%%%%%%
%%%%%%%
\section{A seven cubes problem.}
\label{7cubes}
%%%%%%
%%%%%%
Watson's proof in \cite{Wa} relies on the fact that, for $X>\exp(q^{1/100})$, the existence of a prime $p\equiv a (\text{mod}\,q)$ in the interval $[X,1.01\,X]$ makes it possible to write a sufficiently large integer $n$ as a sum of seven cubes. And the size of the smallest of these $n$'s depends on the size of $X$.
We will follow the latest version of this algorithm, due to Ramar\'e (\cite{R}).
%%%%%%
\subsection{A modified form of Watson's lemma (Lemma $5$ of \cite{MC})}\label{sectionWatson}
The next lemma provides conditions for an integer to be a sum of seven cubes.
\begin{lemma}[Lemma $2.1$ of \cite{R}]
Let $n$, $a$, $u$, $v$ and $w$ be positive integers and $t$ a non negative integer. We assume that
\begin{eqnarray}
&& \label{(1)}
1 \le u \le v \le w \le (3/4)^{1/3} uv/24,\\
&& \label{(2)}
gcd(uvw,6n) = 1 \hbox{ and } a \hbox{ is odd},\\
&& \label{(3)}
u, v, w \hbox{ and } a \hbox{ are pairwise co-prime}, \\
&& \label{(4)}
n-t^3 \equiv 1\ [2] \\
&& \label{(5)}
n-t^3 \equiv 0\ [3a] \\
&& \label{(6)} \left\lbrace
\begin{array}{l}
4(n-t^3) \equiv v^6w^6a^3\ [u^2] \\
4(n-t^3) \equiv u^6w^6a^3\ [v^2] \\
4(n-t^3) \equiv u^6v^6a^3\ [w^2]
\end{array}
\right.
\end{eqnarray}
Set $\d = \left( 1+(w/u)^6 +(w/v)^6 \right)/4$. If
\begin{equation}\label{(7)}
0 \le \frac{uv}{6w} \left( \frac{n}{u^6v^6a^3} - \d - \frac34 \right)^{1/3} \le \frac{t}{6uvwa} \le \frac{uv}{6w} \left( \frac{n}{u^6v^6a^3} - \d \right)^{1/3}
\end{equation}
then $n$ is a sum of seven non-negative cubes.
\end{lemma}
%%%%%%%
%%%%%%%
\subsection{Reducing to finding a prime in an arithmetic progession}\label{reducing}
Suppose the integer $n$ is given.
We need to find $u$, $v$, $w$, $a$ and $t$ such that the conditions of our lemma are fulfilled.\\
Let $u,v,w$ be prime numbers $\equiv 5\ [6]$ that are satifying \myref{(1)} and are coprime with $n$.
Then $(4n)/(v^6w^6)$ is a cube, modulo $u^2$. We have the same for $(4n)/(u^6w^6)$ modulo $v^2$ and $(4n)/(u^6v^6)$ modulo $w^2$.
This is easy to prove, knowing that, if $p$ is prime $\equiv 5\ [6]$, then every invertible residue class modulo $p$ is a cube modulo $p^2$.\\
Moreover $u^2$, $v^2$ and $w^2$ are pair-wise coprime and, by the Chinese remainder theorem, there exists an integer $a^{\prime}$ such that
\begin{equation}
\label{inter}
\left\lbrace \begin{array}{l}
4n \equiv (a' v^2 w^2)^3\ [u^2] \\
4n \equiv (a' u^2 w^2)^3\ [v^2] \\
4n \equiv (a' u^2 v^2)^3\ [w^2]
\end{array}\right.
\end{equation}
We choose $a$ to be
\(a \equiv a' \ [u^2v^2w^2]\),
so that we can replace $a'$ by $a$ in the system \myref{inter}.
Also we can choose $a$ to be prime and
\(a \equiv 5\ [6]\),
so that we are insured that there exists an integer $n$ cubic
modulo $3a$. 
We deduce that
\begin{cond}\label{C1}
There exists a prime $a$ such that $a \equiv \ell \  [6u^2v^2w^2]$.
\end{cond}
Since the integers $u$, $v$, $w$ and $6a$ are coprime, there
exists integers $t$ satisfying:
\begin{equation*}
t^3 \equiv n \ [3a] ,\ t^3 \equiv n-1\ [2] ,\ t\equiv 0\ [uvw].
\end{equation*}
Up to now, the conditions \myref{(1)} to \myref{(6)} are
satisfied. In order to find $\frac{t}{6auvw}$ bounded as in
\myref{(7)}, we need to add some conditions on $a$, namely that
\begin{cond}
\label{C2}
$\dfrac{Y}{\k}  \le a \le  Y$, 
\begin{equation*}
\hbox{where}\ 
Y:= \frac{n^{1/3}}{u^2v^2(3/4 + \d)^{1/3}}
,\k^3 := \dfrac1{3/4 + \d} \left[ \left( \frac{uv}{24w(\rho +1)}\right) ^{3/2} + \d\right] 
,\rho:=\frac{1}{6uvwa}.
\end{equation*}
\end{cond} 
More explanation is provided on pp $377$-$378$ of \cite{R}. We
will see that Theorem \ref{mainthm} insures us of the existence of
a prime $a$ satisfying conditions \ref{C1} and \ref{C2}. However,
this theorem is established for non-exceptional moduli. We explain
how to avoid the case of exceptional zeros in the next section.
%%%%%%%
%%%%%%%
\subsection{Creating a non-exceptional modulus}
\begin{theorem}[Theorem $2$ of \cite{RR}]\label{ramrum}
For all $q \le 72$, and all $a$ prime to $q$, uniformly for $1 \le x \le 10^{10}$,
\begin{equation*}
\max_{1\le y \le x} \left| \theta(y;q,a) - \frac{y}{\phi(q)}\right| \le 2.072 \sqrt{x}.
\end{equation*}
\end{theorem}
\begin{lemma}
There are more than $12$ prime numbers coprime to $n$ and congruent to $5$ modulo $6$ lying in the interval $[0.521\log n, 2.562 \log n]$ if $\log n $ is larger than $68\,509$.
\end{lemma}
\begin{proof}
See proof of Lemma $4.5$ of \cite{R}.
The constants $c_1$ and $c_2$ are choosen to optimize the lower bound of $\log n$ given in the equation \myref{N0cubes1} below under the conditions that $c_2-c_1>\phi(6)$ and
\begin{equation*}
\left(\frac{c_2-c_1}{2}-1\right) \log n - 2.072 \left(\sqrt{c_1}+\sqrt{c_2}\right)\sqrt{\log n} \ge 
12 \log \left(c_2 \log n\right) .
\end{equation*}
\end{proof}
We note $c_3=\left(\frac{c_2}{c_1}\right)^{\frac1{2}}$. Then we deduce by the pigeon hole principle that
there exists an interval
$[A,c_3 A]$ with $A$ in $ [c_1\log n, c_2/c_3 \log n ]$ which contains more than $6$ primes co-prime to $n$ and
congruent to $5$ modulo $6$. We denote $u_1<v_1<w_1<u_2<v_2<w_2$, $6$ of these primes, $k_1=3(u_1v_1w_1)^2$ and $k_2=3(u_2v_2w_2)^2$. To prove that one of the coprime integers $k_1$ and $k_2$ has to be non-exceptional, we use Theorem \ref{repulsion} and the inequalities
\begin{equation*}
k_1^{2.12} \ge \left(3\left(c_1\log n \right)^6\right)^{2.12}
> 3\left(c_2\log n\right)^6
\ge k_2
\end{equation*}
for $n\ge 150$. 
We denote simply $k$ the non-exceptional modulus and $u$, $v$, $w$
the associated integers in $[A,c_3 A]$. It remains to find for
which $n$, the interval $[Y/\k,Y]$ satisfies the hypothesis of
Theorem \ref{mainthm}.
%%%%%%%
%%%%%%%
\subsection{Finding a prime in a progression with a large modulus}
\begin{lemma}
Assume $\log n \ge 71\,000$. For any invertible residue
class $l$ modulo $k$, there is a prime $a$, congruent to $l$ modulo $k$
contained in $[Y/\k,Y]$.
\end{lemma}
\begin{proof}
We use the bounds:
\begin{equation*} c_1\log n \le A \le u,v,w \le c_3 A \le c_2 \log n ,\ 
\frac14+\frac{1}{2c_3^6} \le \d \le \frac14+\frac{c_3^6}2 ,\ 
 \rho \le \frac1{6(c_1\log n)^3}.
\end{equation*} 
We deduce that $\displaystyle{\k \ge \k_0(n)}$ and $\displaystyle{Y\ge Y_0(n)}$, with
\begin{align*}
& \k_0(n)^3 
:= \frac1{1+\frac{c_3^6}{2}} \left( \left( \frac{c_1\log n}{24c_3\left( \frac1{6(c_1\log n)^3} +1\right)}\right) ^{3/2} + \frac14+\frac{1}{2c_3^6}  \right),
\\
& Y_0(n) 
:= \frac{n^{1/3}}{(c_2\log n)^4 \left(1+\frac{c_3^6}{2}\right)^{1/3}}.
\end{align*}
%%%%%%%%%%%%%
For the values $k\ge 10^{32}$, $\a=4.3060$ and $\ep=1.9$, the inequality
$Y \ge e^{\a \log^2 k +\ep}$
is satisfied when
\begin{equation}\label{N0cubes1}
\frac{\log n}3  - 4\log(c_2\log n) - \frac13 \log \left(1+\frac{c_3^6}2 \right)\ge \a \log^2 \left(3 (c_2 \log n)^6 \right) + \ep,
\end{equation}
that is to say for $\log n\ge 70\,341$. This also warrants $\k_0(n) \ge e^{\ep}$.
\end{proof}
%%%%%%%
%%%%%%%
\begin{tiny}
\begin{table}[htbp]
\caption{} \label{computations} 
\renewcommand\arraystretch{0.95}
\noindent
\begin{equation*}
\begin{array}[t]{@{}c@{}@{}c@{}}
\begin{array}[t]{|c|c|c|c|r|r|}
\hline
q_0 & \ep & u & m  & H  & \a  
\\ \hline
\multirow{6}{*}{$5 . 10^4$} 
& 0.0001 & 0.086 & 14 & 514\,998 & 19.228
\\& 0.001& 0.092 & 13 & 47\,292 & 15.550
\\& 0.01 & 0.098 & 12 & 4\,311 &12.245 
\\& 0.1  & 0.004 & 14 & 528 & 9.4357
\\& 1    & 0.01 & 15 & 57.8 & 6.9684
\\& 10   & 0.037 & 11 & 4.4219 & 4.8430
\\ \hline
\multirow{6}{*}{$10^{10}$} 
& 0.0001 & 0.056 & 20 & 741\,876 & 9.8356
\\& 0.001& 0.058 & 19 & 70\,330 & 8.5912
\\& 0.01 & 0.060 & 18 & 6\,632 & 7.4254
\\& 0.1  & 0.061 & 17 & 630 & 6.3398
\\& 1    & 0.057 & 16 & 62.5 & 5.3418
\\& 10   & 0.028 & 17 & 6.75 & 4.4761
\\ \hline
\multirow{6}{*}{$10^{15}$} 
& 0.0001 &0.043 & 25 & 948\,594 & 7.6121
\\& 0.001&0.045 & 24 & 90\,920 & 6.8799
\\& 0.01 &0.046 & 23 & 8\,713 & 6.1816
\\& 0.1  &0.046 & 22 & 839 & 5.5174
\\& 1    &0.043 & 21 & 83.6 & 4.8905
\\& 10   &0.024 & 22 & 8.84 & 4.3256
\\ \hline
\multirow{6}{*}{$10^{20}$} 
& 0.0001 &0.035 & 30 & 1\,152\,223  & 6.5919
\\& 0.001&0.036 & 29 & 111\,390 & 6.0799
\\& 0.01 &0.037 & 28 & 10\,762 & 5.5864
\\& 0.1  &0.037 & 27 & 1\,045 & 5.1114
\\& 1    &0.035 & 26 & 105 & 4.6565
\\& 10   &0.021 & 27 & 10.9 & 4.2373
\\ \hline
\multirow{6}{*}{$10^{25}$} 
& 0.0001 &0.030 & 35 & 1\,353\,117 & 6.0079
\\& 0.001&0.030 & 34 & 131\,618 & 5.6164
\\& 0.01 &0.031 & 33 & 12\,786 & 5.2364
\\& 0.1  &0.031 & 32 & 1\,247 & 4.8678
\\& 1    &0.029 & 31 & 125 & 4.5116
\\& 10   &0.019 & 32 & 12.9 & 4.1783
\\ \hline
\multirow{6}{*}{$10^{30}$} 
& 0.0001 &0.026 & 40 & 1\,553\,007 & 5.6298
\\& 0.001&0.026 & 39 & 151\,626 & 5.3137
\\& 0.01 &0.026 & 38 & 14\,802 & 5.0053
\\& 0.1  &0.026 & 37 & 1\,449 & 4.7046
\\& 1    &0.025 & 36 & 145 & 4.4123
\\& 10   &0.017 & 37 & 15 & 4.1357
\\ \hline
\multirow{6}{*}{$10^{35}$} 
& 0.0001 &0.023 & 45 & 1\,751\,630 & 5.3649
\\& 0.001&0.023 & 44 & 171\,503 & 5.1002
\\& 0.01 &0.023 & 43 & 16\,791 & 4.8411
\\& 0.1  &0.023 & 42 & 1\,648 & 4.5875
\\& 1    &0.022 & 41 & 165 & 4.3396
\\& 10   &0.016 & 42 & 16.9 & 4.1032
\\ \hline
\multirow{6}{*}{$10^{40}$} 
& 0.0001 &0.020 & 50 & 1\,950\,568 & 5.1688
\\& 0.001&0.021 & 49 & 191\,213 & 4.9414
\\& 0.01 &0.021 & 48 & 18\,763 & 4.7181
\\& 0.1  &0.021 & 47 & 1\,846 & 4.4989
\\& 1    &0.020 & 46 & 185 & 4.2839
\\& 10   &0.014 & 47 & 18.9 & 4.0776
\\ \hline
\multirow{6}{*}{$10^{45}$} 
& 0.0001 &0.018 & 55 & 2\,114\,784 & 5.0178
\\& 0.001&0.019 & 54 & 210\,928 & 4.8185
\\& 0.01 &0.019 & 53 & 20\,736 & 4.6225
\\& 0.1  &0.019 & 52 & 2\,043 & 4.4295
\\& 1    &0.018 & 51 & 204 & 4.2398
\\& 10   &0.013 & 52 & 20.9 & 4.0567 
\\ \hline
\multirow{6}{*}{$10^{50}$} 
& 0.0001 &0.017 & 60 & 2\,342\,931 & 4.8979
\\& 0.001&0.017 & 59 & 230\,659 & 4.7206
\\& 0.01 &0.017 & 58 & 22\,710 & 4.5459
\\& 0.1  &0.017 & 57 & 2\,240 & 4.3737
\\& 1    &0.016 & 56 & 224 & 4.2039
\\& 10   &0.012 & 57 & 22.9 & 4.0394
\\ \hline
\end{array}
\end{array}
\begin{array}[t]{ll}
\begin{array}[t]{|c|c|c|c|r|r|} 
\hline
q_0 & \ep & u & m  & H  & \a  
\\ \hline
\multirow{6}{*}{$10^{55}$} 
& 0.0001 &0.016 & 64 & 2\,538\,632 & 4.8003
\\& 0.001&0.016 & 64 & 250\,167 & 4.6407
\\& 0.01 &0.016 & 63 & 24\,661 & 4.4832
\\& 0.1  &0.015 & 62 & 2\,338 & 4.3276
\\& 1    &0.015 & 61 & 244 & 4.1740
\\& 10   &0.012 & 62 & 24.8 & 4.0247
\\ \hline
\multirow{6}{*}{$10^{60}$} 
& 0.0001 &0.015 & 69 & 2\,731\,576 & 4.7192
\\& 0.001&0.015 & 68 & 269\,639 & 4.5742
\\& 0.01 &0.014 & 68 & 26\,639 & 4.4308
\\& 0.1  &0.014 & 67 & 2\,633 & 4.2890
\\& 1    &0.014 & 66 & 263 & 4.1488
\\& 10   &0.011 & 67 & 26.8 & 4.0121
\\ \hline
\multirow{6}{*}{$10^{65}$} 
& 0.0001 &0.014 & 74 & 2\,927\,544 & 4.6509
\\& 0.001&0.014 & 73 & 289\,140 & 4.5179
\\& 0.01 &0.014 & 72 & 28\,660 & 4.3864
\\& 0.1  &0.013 & 72 & 2\,829 & 4.2562
\\& 1    &0.013 & 71 & 283 & 4.1272
\\& 10   &0.010 & 72 & 28.7 & 4.0011
\\ \hline
\multirow{6}{*}{$10^{70}$} 
& 0.0001 &0.013 & 79 & 3\,122\,581 & 4.5924
\\& 0.001& 0.013 & 78 & 308\,647 & 4.4697
\\& 0.01 & 0.013 & 77 & 30\,511 & 4.3482
\\& 0.1  &0.013 & 76 & 3\,021 & 4.2278
\\& 1    &0.012 & 76 & 302 & 4.1084
\\& 10   &0.010 & 77 & 30.7 & 3.9915
\\ \hline
\multirow{6}{*}{$10^{75}$} 
& 0.0001 &0.012 & 84 & 3\,317\,736 & 4.5418
\\& 0.001&0.012 & 83 & 328\,165 & 4.4280
\\& 0.01 &0.012 & 82 & 32\,465 & 4.3151
\\& 0.1  &0.012 & 81 & 3\,216 & 4.2031
\\& 1    &0.011 & 81 & 322 & 4.0920
\\& 10   &0.009 & 81 & 32.6 & 3.9829
\\ \hline
\multirow{6}{*}{$10^{80}$} 
& 0.0001 &0.011 & 89 & 3\,513\,060 & 4.4976
\\& 0.001&0.011 & 88 & 347\,700 & 4.3914
\\& 0.01 &0.011 & 87 & 34\,417 & 4.2860
\\& 0.1  &0.011 & 86 & 3\,311 & 4.1814
\\& 1    &0.011 & 86 & 341 & 4.0774
\\& 10   &0.009 & 86 & 34.5 & 3.9753
\\ \hline
\multirow{6}{*}{$10^{85}$} 
& 0.0001 &0.011 & 94 & 3\,704\,920 & 4.4587
\\& 0.001&0.011 & 93 & 366\,878 & 4.3591
\\& 0.01 &0.011 & 92 & 36\,335 & 4.2603
\\& 0.1  &0.010 & 91 & 3\,607 & 4.1621
\\& 1    &0.010 & 90 & 361 & 4.0645
\\& 10   &0.008 & 91 & 35.5 & 3.9684
\\ \hline
\multirow{6}{*}{$10^{90}$} 
& 0.0001 &0.010 & 98 & 3\,900\,103 & 4.4240
\\& 0.001&0.010 & 98 & 386\,427 & 4.3331
\\& 0.01 &0.010 & 97 & 38\,290 & 4.2373
\\& 0.1  &0.010 & 96 & 3\,799 & 4.1448
\\& 1    &0.010 & 95 & 380 & 4.0528
\\& 10   &0.008 & 96 & 37.4 & 3.9623 
\\ \hline
\multirow{6}{*}{$10^{95}$} 
& 0.0001 &0.010 & 103 & 4\,091\,636  & 4.3931
\\& 0.001&0.010 & 102 & 405\,566 & 4.3046
\\& 0.01 &0.010 & 101 & 40\,204 & 4.2168
\\& 0.1  &0.009 & 101 & 3\,995 & 4.1293
\\& 1    &0.009 & 100 & 399 & 4.0423
\\& 10   &0.008 & 101 & 40.3 & 3.9565 
\\ \hline
\multirow{6}{*}{$10^{100}$} 
& 0.0001 &0.009 & 108 & 4\,287\,331  & 4.3652
\\& 0.001&0.009 & 107 & 425\,137 & 4.2815
\\& 0.01 &0.009 & 106 & 41\,161 & 4.1982
\\& 0.1  &0.009 & 105 & 4\,186 & 4.1153
\\& 1    &0.009 & 105 & 419 & 4.0328
\\& 10   &0.007 & 106 & 42.3 & 3.9513
\\ \hline
\end{array}
\end{array}
\end{equation*}
\end{table}
\end{tiny}
%%%%%%%
%%%%%%%
\bibliographystyle{amsplain}

\end{document}